\documentclass{amsart}

\usepackage{amsmath}
\usepackage{amsfonts}
\usepackage[utf8]{inputenc}
\usepackage[numbers]{natbib}
\usepackage{graphicx}
\usepackage{mathtools}
\usepackage{amssymb}
\usepackage{amsthm}
\usepackage{tikz-cd} 
\usepackage{enumerate}
\usepackage[english]{babel}
\usepackage{hyperref}
\usepackage{tikz}
\usetikzlibrary{positioning} 
\usepackage{nomencl}

\theoremstyle{definition}
\newtheorem{cor}{Corollary}[section] 
\newtheorem{prop}[cor]{Proposition}
\newtheorem{theorem}[cor]{Theorem}

\newtheorem{conj}[cor]{Conjecture}
\newtheorem{defi}[cor]{Definition}
\newtheorem{bsp}[cor]{Example}
\newtheorem{lemma}[cor]{Lemma}

\numberwithin{equation}{section}

\DeclareMathOperator{\simk}{sc}

\DeclareMathOperator{\aut}{Aut}

\DeclareMathOperator{\pgl}{PGL}

\DeclareMathOperator{\irr}{Irr}

\DeclareMathOperator{\mor}{Mor}

\DeclareMathOperator{\iso}{Iso}

\DeclareMathOperator{\syl}{Syl}

\DeclareMathOperator{\inj}{Inj}

\DeclareMathOperator{\ad}{ad}
\DeclareMathOperator{\res}{Res}

\DeclareMathOperator{\trr}{Tr}

\DeclareMathOperator{\homo}{Hom}

\DeclareMathOperator{\br}{Br}

\begin{document}


\title{Reduction Theorems for generalised block fusion systems}


\author{Patrick Serwene}
\address{Technische Universität Dresden, Faculty of Mathematics, 01062 Dresden, Germany}
\email{patrick.serwene@tu-dresden.de}





\begin{abstract}
We generalise Brauer's Third Main Theorem, the First and Second Fong Reduction to generalised block fusion systems and apply the Second Fong Reduction to extend a result by Cabanes about the non-exoticity of fusion systems of unipotent blocks for finite groups of Lie type.
\end{abstract}






\maketitle

\section{Introduction}
Let $p$ be a prime number, $k$ an algebraically closed field of characteristic $p$, $P$ a finite $p$-group and $\mathcal F$ a fusion system on $P$. By fusion system we always mean saturated fusion system. Recall that $\mathcal F$ is said to be realisable if $\mathcal F=\mathcal F_P(G)$ for a finite group $G$ and $P \in \syl_p(G)$, otherwise $\mathcal F$ is said to be exotic. Also recall that if $\mathcal F=\mathcal F_{(P,e_P)}(G,b)$ for a finite group $G$ having a $p$-block $b$ with maximal $b$-Brauer pair $(P,e_P)$, $\mathcal F$ is said to be block-realisable, otherwise $\mathcal F$ is said to be block-exotic.\\
The following fact is a consequence of Brauer's Third Main Theorem (see \cite[Theorem 3.6]{radha}): If $G$ is a finite group and $b$ is the principal $p$-block of $kG$, i.e. the block corresponding to the trivial character, with maximal $b$-Brauer pair $(P,e_P)$, then $P \in \syl_p(G)$ and $\mathcal{F}_{(P,e_P)}(G,b)=\mathcal{F}_P(G)$. In particular, any realisable fusion system is block-realisable. The converse is still an open problem and has been around in the form of the following conjecture for  a  while, see \cite[Part IV,7.1]{ako} and \cite[9.4]{david}:
\let\thefootnote\relax\footnotetext{This work formed part of the author’s PhD thesis at City, University of London under the
supervision of Prof. Radha Kessar.}

\begin{conj}
\label{XX}
\textit{If $\mathcal F$ is an exotic fusion system, then $\mathcal F$ is block-exotic.}
\end{conj}

When studying block fusion systems, and proving theorems about the aforementioned conjecture, one needs a category which generalises block fusion systems, which we call generalised block fusion systems. The main reason is that block fusion systems do not work well with respect to Clifford Theory. Let us first recall an important term.

\begin{defi}
\textit{Let $G$ be a finite group and $N \unlhd G$. Let $b$ be a block of $kG$ and $c$ be a block of $kN$. We say that $b$ \textbf{covers} $c$ if $bc \neq 0$.}
\end{defi} 

Let $G$ be a finite group with $N \unlhd G$ and $c$ a block of $kN$ covered by a block $b$ of $kG$. Let $P$ be a defect group for $b$, then $P \cap N$ is a defect group for $c$, but the fusion system of $c$ is not necessarily contained in the fusion system of $b$, which means that block fusion systems do in general not allow for a descent to normal subgroups. However, one can circumvent this problem by using generalised block fusion systems. These were introduced in \cite{ks} and play a key role in the main result of \cite{ks}, as well as in the Reduction Theorems of \cite{serwene}. They are defined as follows. See Section \ref{general} for details.

\begin{defi}
\textit{Let $H$ be a finite group, $G \unlhd H$ and $c$ be a $H$-stable block of $kG$.  Let $(P,e_P)$ be a maximal $(c,H)$-Brauer pair. Let $(Q,e_Q) \leq (P,e_P)$. Denote by $\mathcal F_{(P,e_P)}(H,G,c)$ the category on $P$ with morphisms consisting of all injective group homomorphisms $\varphi: Q \rightarrow R$ for which there is some $g \in H$ such that $\varphi(x)$ = $^gx$ for all $x \in Q$ and $^g(Q,e_Q) \leq (R,e_R)$. We refer to this category as generalised block fusion systems}.
\end{defi}

Generalised block fusion systems capture the relationship between covered blocks of normal subgroups and the covering blocks of their overgroups. We explain these relationships in detail in the background chapter. We make the following conjecture.

\begin{conj}
\label{XZ}
\textit{Let $H$ be a finite group with normal subgroup $G$ having an $H$-stable $p$-block $b$ with maximal $(b,H)$-Brauer pair $(P,e_P)$. Then the generalised block fusion system $\mathcal F_{(P,e_P)}(H,G,b)$ is non-exotic.}
\end{conj}

When $H=G$, then $\mathcal F_{(P,e_P)}(H,G,b)$ is the usual block fusion system of $b$, so $\mathcal F_{(P,e_P)}(H,G,b)=\mathcal F_{(P,e_P)}(G,b)$. In particular, Conjecture \ref{XZ} implies Conjecture \ref{XX}. We prove three key results for block fusion systems for this more general category, namely Brauer's Third Main Theorem and the First and Second Fong Reduction. The generalised version of Brauer's Third Main Theorem is the following result. If $G=H$, we get the result for block fusion systems.

\begin{theorem}
\label{a}
\textit{Let $G \unlhd H$, $b$ the principal block of $kG$, $(P,e_P)$ a maximal $(b,H)$-Brauer pair, then $P \in \syl_p(H)$, $e_P$ is the principal block of $kC_G(P)$ and
 $ \mathcal F_P(H)=\mathcal F_{(P,e_P)}(H,G,b)$.}
\end{theorem}

The next result is the Generalised First Fong Reduction. If $G=H$ in the following theorem, we obtain the original result for block fusion systems. This special case was proved in \cite[Part IV, Proposition 6.3]{ako}. Let $G \unlhd H$, $c$ be a block of $kG$ and $I_H(c)=\{h \in H \mid ^hc=c \}$, the inertia subgroup of $c$ in $H$. If $b$ is a block of $kH$ covering $c$, then there is a block of $kI_H(c)$ corresponding to $b$ which we call the Fong correspondent of $b$.

\begin{theorem}
\label{b}
\textit{Let $\mathcal F$ be a fusion system on a $p$-group $P$ and let $G$, $H$ be finite groups such that $G \unlhd H$. Let $b$ be an $H$-stable block of $kG$ with $\mathcal F=\mathcal F_{(P,e_P)}(H,G,b)$. Let $N$ be a normal subgroup of $H$ contained in $G$ and $c$ be a block of $kN$ which is covered by $b$. Then $\mathcal F=\mathcal F_{(P,\widetilde{e_P})}(I_H(c),I_G(c),\tilde b)$, where $\tilde b$ is the Fong correspondent of $b$ in $I_G(c)$ and $(P,\widetilde{e_P})$ is a maximal $(\tilde b, I_H(c))$-Brauer pair}.
\end{theorem}

Finally, we state the Generalised Second Fong Reduction. If $G=H$, we obtain the original theorem for block fusion systems. This special case was proved in  \cite[Part IV, Theorem 6.6]{ako}.

\begin{theorem}
\label{c}
\textit{Let $M \leq H$ such that $|H:M|_{p}=1$ and let $A$ be a normal subgroup of $H$ contained in $M$. Let $c$ be an $H$-stable block of $kA$ of defect zero and $d$ be an $H$-stable block of $kM$ covering $c$ with maximal $(d,H)$-Brauer pair $(P,e_P)$. Then there exists a $p'$-central extension $\widetilde H$ of $H/A$ and a block $\tilde d$ of $\widetilde M$ with maximal $(\tilde d, \widetilde H)$-Brauer pair $(\widetilde P,e_P')$, with $\widetilde P \cong P$, where $\widetilde M$ is the full inverse image of $M/A$ in $\widetilde H$ such that $\tilde d$ is $\widetilde H$-stable and $\mathcal F_{(P,e_P)}(H,M,d) \cong \mathcal F_{(\widetilde P,e_P')}(\widetilde H,\widetilde M,\tilde d)$}.
\end{theorem}

Concerning Conjectures \ref{XX} and \ref{XZ}, a first goal would be to prove these for all finite (quasi-)simple groups. We are interested in the case of finite simple groups of Lie type. If the natural characteristic of the group and $p$ coincide, Conjecture \ref{XX} is known to hold, see \cite[Theorem 6.18]{caen}. If $p$ is not the natural characteristic, Cabanes proved Conjecture \ref{XX} via constructive argument, i.e. providing a group which realises the fusion system, in \cite{oh} for unipotent $p$-blocks, where $p$ is at least $7$. The next main result extends the result of Cabanes to generalised block fusion systems.

\begin{theorem}
\label{kombucha}
\textit{Let $\mathbf X$ be a connected, reductive group with associated Frobenius $F$, $\mathbf G \leq \mathbf X$ an $F$-stable Levi subgroup having a unipotent block $b$ and let $\mathbf N \leq N_{\mathbf X}(\mathbf G)$ be $F$-stable. We furthermore assume $p \geq 7$, $|\mathbf N^F/\mathbf G^F|_p=1$ and that $b$ is $\mathbf N^F$-stable. Let $(P,e_P)$ be a maximal $(b, \mathbf N^F)$-Brauer pair. The generalised block fusion system $\mathcal F_{(P,e_P)}(\mathbf N^F,\mathbf G^F,b)$ is non-exotic}.
\end{theorem}

The case $\mathbf G=\mathbf N$ is Cabanes' Theorem. Let $(\mathbf X,F)$ be in duality with $(\mathbf X^\ast,F)$. If $d$ is a block of $\mathbf X^F$, then there exists a semisimple $p'$-element $s$ in $\mathbf {X^\ast}^F$ such that $\irr(\mathbf X^F,d) \cap \irr(\mathbf X^F,s) \neq 0$. Let $\mathbf G^\ast$ be a Levi subgroup of $\mathbf X^\ast$ minimal with respect to the property of containing $C^\circ_{\mathbf X^\ast}(s)$ and denote by $\bf G$ its dual. In \cite{bdr}, Bonnafé, Dat and Rouquier prove that, in many cases, the fusion system of $d$ is equivalent to the fusion system of a block $c$ of a subgroup of $\mathbf X^F$ containing $\mathbf G^F$ as a normal subgroup, where $c$ covers a unipotent block of $\mathbf G^F$. We can use Theorem \ref{kombucha} to prove that the generalised block fusion system related to this normal inclusion is non-exotic and hope to use this in further work to tackle Conjectures \ref{XX} and \ref{XZ}. To be specific, a next step would be to find a normal pair $G' \unlhd H'$ of finite groups determined by $(\mathbf X^F,d)$ and some $H'$-stable unipotent block $b'$ of $G'$ such that generalised block fusion system of $(H',G',b')$ coincides with the fusion system of $(\mathbf X^F,d)$. \\
In the next section we recall some definitions and results about fusion systems, in particular also block fusion systems and generalised block fusion systems. In Section \ref{drei}, we first state the original versions of Theorems \ref{a}, \ref{b} and \ref{c}, and then prove the new versions. Finally, in Section \ref{funf} we apply Theorem \ref{b} to prove Theorem \ref{kombucha}.

\section{Notation and Terminology}
\subsection{Introduction to Fusion Systems}
In this section, we recall the definition and some key properties of fusion systems. Let $p$ be a prime number. Note that by $p$-group we always mean finite $p$-group. For two groups $Q, R$ contained in a group $G$, we write $\homo_G(Q,R)$ for all group homomorphisms between $Q$ and $R$ which are induced by conjugation with an element in $G$.

\begin{defi}
\label{jabloka}
\textit{$(a)$ Let $p$ be a prime and $P$ be a $p$-group. A \textbf{fusion system} is a category $\mathcal F$ on $P$, such that for all $Q, R \leq P$ we have:\\
$(i) \ \homo_P(Q,R) \subseteq \homo_{\mathcal F}(Q,R) \subseteq \inj(Q,R),$\\
$(ii)$ each homomorphism in $\mathcal F$ is the composition of an $\mathcal F$-isomorphism and an inclusion.\\
$(b)$ Let $\mathcal F$ be a fusion system on a $p$-group $P$. Two subgroups $Q,R \leq P$ are called \textbf{$\mathcal F$-conjugate} if they are isomorphic as objects of the category $\mathcal F.$\\
$(c)$ A subgroup $Q \leq P$ is called \textbf{fully automised} in $\mathcal F$ if $\aut_P(Q) \in \syl_p(\aut_{\mathcal F}(Q)).$\\
$(d)$ A subgroup $Q \leq P$ is called \textbf{receptive} in $\mathcal F$ if for each $R \leq P$ and each $\varphi \in \iso_{\mathcal{F}}(R,Q)$, $\varphi$ has an extension to the group $N_{\varphi}:=N_{\varphi}^{\mathcal F}:=\{ g \in N_P(R) \mid$ $^{\varphi} c_{g} \in \aut_P(Q) \}$.\\
$(e)$ A fusion system is called \textbf{saturated} if each subgroup of $P$ is $\mathcal F$-conjugate to a subgroup which is fully automised and receptive in $\mathcal F$.}
\end{defi}

In many applications, it is crucial for fusion systems to be saturated, since fusion systems only satisfying part $(a)$ of the previous definition are too general. For convenience, we drop the term saturated, and mean saturated fusion system whenever we say fusion system. In the literature, fusion system means categories satisfying only part $(a)$ from Definition \ref{jabloka}.

\begin{theorem}
\cite[Theorem 2.11]{markus}
\textit{Let $G$ be a finite group with $P \in \syl_p(G)$. We denote the category on $P$ with morphisms consisting of homomorphisms induced by conjugation by elements in $G$ by $\mathcal{F}_P(G)$. Then $\mathcal F_P(G)$ is a fusion system on $P$}.
\end{theorem}

If a fusion system is of the form $\mathcal F_P(G)$ for a finite group $G$ and $P \in \syl_p(G)$, we call it realisable, otherwise we call it exotic. Furthermore, we say that a fusion system on a $p$-group $P$ is trivial if $\mathcal F=\mathcal F_P(P)$.

\begin{defi}
\textit{Let $\mathcal F$ be a fusion system on a $p$-group $P$.\\
$(a)$ If $C_{P}(Q')=Z(Q')$ for each $Q' \leq P$ which is $\mathcal F$-conjugate to $Q$, then $Q$ is called \textbf{$\mathcal F$-centric}. Define $\mathcal F^c$ to be the full subcategory of $\mathcal F$ whose objects are the $\mathcal F$-centric subgroups of $P$.\\
$(b)$ A proper subgroup $H$ of a finite group $G$ is called \textbf{strongly $p$-embedded} if $H$ contains a Sylow $p$-subgroup $P$ of $G$ and $P \neq 1$ but $^x P \cap H=1$ for any $x \in G \setminus H$.\\
$(c)$ A subgroup $Q \leq P$ is called \textbf{fully $\mathcal F$-normalised} if $|N_P(R)| \leq |N_P(Q)|$ for any $R \leq P$ with $R \cong Q$ in $\mathcal F$.\\
$(d)$ We call a subgroup $Q \leq P$ \textbf{$\mathcal F$-essential} if $Q$ is $\mathcal F$-centric and fully normalised in $\mathcal F$, and if $\aut_{\mathcal F}(Q)/\aut_Q(Q)$ has a strongly $p$-embedded subgroup}.
\end{defi}

Note that if $\mathcal F$ is a fusion system on $P$, an $\mathcal F$-essential subgroup of $P$ is always a proper subgroup since $\aut_{\mathcal F}(P)/\aut_P(P)$ is a $p'$-group.\\

Let $n \in \mathbb N_{\geq 1}$ and $i \in \{1, \dots, n\}$. If $\mathcal F$ is a fusion system on $P$ and $M_i \subseteq \aut_{\mathcal F}(Q_i)$ for some $Q_i \leq P$, we denote by $\langle M_1, \dots, M_n \rangle$ the smallest (not necessarily saturated) subsystem of $\mathcal F$ on $P$ such that its morphisms contain all the sets $M_i$ for $i=1,\dots,n$. The following theorem tells us that the structure of a fusion system $\mathcal F$ on $P$ is determined by the automorphisms of $\mathcal F$-essential subgroups of $P$ and $P$ itself.

\begin{theorem}
\label{alperin} (Alperin's Fusion Theorem)
\cite[Part I, Theorem 3.5]{ako}
\textit{Let $\mathcal F$ be a fusion system on a $p$-group $P$. Then $\mathcal F=\langle \aut_{\mathcal F}(Q) \mid Q=P$ or $Q$ is $\mathcal F$-essential $\rangle$.}
\end{theorem}

We can define a substructure similar to normal subgroups for fusion systems.

\begin{defi}
\textit{
Let $\mathcal F$ be a fusion system on a $p$-group $P$ and $\mathcal E \subseteq \mathcal F$ be a subcategory of $\mathcal F$ which is a fusion system itself on some subgroup $P' \leq P$.\\
$(a)$ A subgroup $Q \leq P$ is called \textbf{strongly $\mathcal F$-closed}, if $\varphi(R) \subseteq Q$ for each $\varphi \in \homo_{\mathcal F}(R,P)$ and each $R \leq Q$.\\
$(b)$ If $P'$ is normal in $P$ and strongly $\mathcal F$-closed, $^\alpha \mathcal E=\mathcal E$ for each $\alpha \in \aut_{\mathcal F}(P')$ and for each $Q \leq P'$ and $\varphi \in \homo_{\mathcal F}(Q,P')$, there are $\alpha \in \aut_{\mathcal F}(P')$ and $\varphi_0 \in \homo_{\mathcal E}(Q, P')$ with $\varphi = \alpha \circ \varphi_0$, then $\mathcal E$ is called \textbf{weakly normal} in $\mathcal F$, denoted $\mathcal E \dot{\unlhd} \mathcal F$.\\
$(c)$ If $\mathcal E$ is weakly normal and in addition, we have that each $\alpha \in \aut_{\mathcal E}(P')$ has an extension $\overline{\alpha} \in \aut_{\mathcal F}(P'C_P(P'))$ with $[\overline{\alpha}, C_P(P')] \leq Z(P')$, then we call $\mathcal E$ \textbf{normal} in $\mathcal F$ and write $\mathcal E \unlhd \mathcal F$}.
\end{defi}

\subsection{Fusion systems of blocks}
In the previous section, we have seen that every finite group realises a fusion system on its Sylow $p$-subgroups. Next, we see that fusion systems can also be induced by block idempotents of finite groups. We refer to such idempotents simply as blocks. Fix $k$ to be an algebraically closed field of characteristic $p$ for the rest of this chapter.

\begin{defi}
\textit{Let $G$ be a finite group and $b$ a block of $kG$. We denote the set of blocks of $kC_G(Q)$ for some $p$-subgroup $Q$ of $G$ by $\mathcal B(Q)$. A \textbf{Brauer pair} is a pair $(Q,f)$ where $Q$ is a $p$-subgroup of $G$ and $f$ is a block of $kC_G(Q)$}.
\end{defi}

Note that $G$ acts on the set of Brauer pairs by conjugation. We recall the Brauer map to see how Brauer pairs form a poset.

\begin{defi}
\textit{Let $G$ be a finite group and $Q \leq G$. For an element $a = \sum\limits_{g \in G} \alpha_g g \in kG$, set $\br^G_Q(a):=\sum\limits_{g \in C_G(Q)} \alpha_g g$}.
\end{defi}

\begin{prop} \cite[Proposition 2.2]{radha}
\textit{Let $G$ be a finite group and let $Q$ be a $p$-subgroup of $G$, henceforth written as $Q \leq_p G$. Then for any $a, a' \in (kG)^Q$, $\br^G_Q(aa')=\br^G_Q(a) \br^G_Q(a')$. Consequently, the map $\br^G_Q: (kG)^Q \rightarrow kC_G(Q), a \mapsto \br^G_Q(a)$ is a surjective homomorphism of $k$-algebras.}
\end{prop}

\begin{defi}
\textit{Let $G$ be a finite group with $N \unlhd G$. We call a block $b$ of $kN$ \textbf{$G$-stable} if $^gb=b$ for all $g \in G$.}
\end{defi}

\begin{defi}
\textit{Let $G$ be a finite group, $Q,R \leq G$ and $(Q,f)$ and $(R,e)$ be Brauer pairs. Then\\
$(a) \ (Q,f) \unlhd (R,e)$ if $Q \unlhd R$, $f$ is $R$-stable and $\br^G_R(f)e=e,$\\
$(b) \ (Q,f) \leq (R,e)$ if $ Q \leq R$ and there exist Brauer pairs $(S_i,d_i), 1 \leq i \leq n$, such that $(Q,f) \unlhd (S_1, d_1) \unlhd (S_2, d_2) \unlhd \dots \unlhd (S_n,d_n) \unlhd (R,e)$}.
\end{defi}

Let $(R,e)$ be a Brauer pair and let $Q \leq R$. The idempotent $f$ such that $(Q,f) \leq (R,e)$ as in the previous definition, is actually uniquely determined:

\begin{theorem} (Alperin--Broué) \cite[Theorem 2.9]{radha}
\label{alperinbro}
\textit{Let $G$ be a finite group, $R \leq G$ and let $(R,e)$ be a Brauer pair. For any $Q \leq R$, there exists a unique $f \in \mathcal B(Q)$ with $(Q,f) \leq (R,e)$. Furthermore, if $Q \unlhd R$, then $f$ is the unique element of $\mathcal B(Q)^R$ with $\br^G_R(f)e=e$. The conjugation action of $G$ on the set of Brauer pairs preserves $\leq$.}
\end{theorem}

\begin{defi}
\textit{Let $G$ be a finite group and $b$ a block of $kG$.\\
$(a)$ A \textbf{$b$-Brauer pair} is a Brauer pair $(R,e)$ such that $(1,b) \leq (R,e)$, or equivalently it is a Brauer pair $(R,e)$ such that $\br^G_R(b)e=e.$\\
$(b)$ We denote the blocks $e$ of $kC_G(R)$ such that $(1,b) \leq (R,e)$ by $\mathcal B(R,b).$\\
$(c)$ A \textbf{defect group} of $b$ is a $p$-subgroup $P$ of $G$ maximal such that $\br^G_P(b) \neq 0$.}
\end{defi}

Note that the group $G$ acts by conjugation on the set of $b$-Brauer pairs. Furthermore, some $p$-subgroup $P$ of $G$ is a defect group of $b$ if and only if there is a maximal pair $(P,e)$ such that $(1,b) \leq (P,e)$. We refer to such a pair as a maximal $b$-Brauer pair.\\

We record some facts about covered blocks which we need in later chapters.

 \begin{lemma}
 \label{rtrt}
 \textit{Let $G$ be a finite group, $S \leq_p G$, $(S,e_S)$ and $(S,f_S)$ be Brauer pairs such that $f_S$ covers $e_S$. Let $e_U$ and $f_U$ be the unique idempotents such that $(U,e_U) \leq (S,e_S)$, respectively $(U,f_U) \leq (S,f_S)$. Suppose that $e_U f_U \neq 0$ for some $U \leq S$. Then $e_Q f_Q \neq 0$ for any $Q \leq U$.}
 \end{lemma}
 
 \textit{Proof}. The second and third paragraph of the proof of \cite[Theorem 3.5]{ks} are spent proving exactly this statement.
 \hfill$\square$
 
 \begin{lemma}  
\label{coverd}
\textit{Let  $G$ be a finite group, $N$ a normal subgroup of $G$, $d$  a  block of $kG$  and $c$ an $G$-stable block of $kN$. Suppose that  there exists  a  $d$-Brauer pair $ (Q, e)  $ such that  $ N \leq C_G(Q)
 $ and $e$ covers  $ c$. Then $d$ covers $c$.}
 \end{lemma}

\textit{Proof}  Since  $e$ covers  $c$ and    $c$ is $C_G(Q)$-stable, $ce = e $, see \cite[Proposition 6.8.2(ii)]{block}. So,
 $dce= de  $.    On the other hand, since $ (Q, e) $ is a  $d$-Brauer pair,  $ \br^G_Q(de)  =  \br^G_Q(d) \br^G_Q (e )  = \br^G_Q(d) e  =  e \ne 0  $. Hence,   $  dce =  de \ne 0  $  which implies that $ dc \neq 0 $. \hfill$\square$
 
  \begin{lemma}
 \label{feit}
 \cite[Chapter V, Lemma 3.5]{feit}
 \textit{Let $G$ be a finite group, $N$ a normal subgroup of $G$ such that $G/N$ is a $p$-group. If $b$ is a block of $kN$, then there is a unique block of $kG$ that covers $b$.}
 \end{lemma}
 
 \begin{defi}
 \textit{Let $b$ be a block of $kG$ and $(P,e_P)$ be a maximal $b$-Brauer pair. For a subgroup $Q \leq P$, denote by $e_Q$ the unique block such that $(Q,e_Q) \leq (P,e_P)$. Denote the category on $P$ whose morphisms consist of all injective group homomorphisms $\varphi: Q \rightarrow R$ for which there is some $g \in G$ such that $\varphi(x)$ = $^gx$ for all $x \in Q$ and $^g(Q,e_Q) \leq (R,e_R)$ by $\mathcal F_{(P,e_P)}(G,b)$}.
 \end{defi}

\begin{theorem}
\label{domingo}
\cite[Theorem 3.9(i)]{radha} \textit{Keep the notation of the previous definition. The category $\mathcal F_{(P,e_P)}(G,b)$ is a fusion system on $P$}.
\end{theorem}

If a fusion system is of the form as in the previous theorem, we call it block-realisable, otherwise we call it block-exotic. Furthermore, if a fusion system $\mathcal F$ is of this form, we call the respective block $b$ $\mathcal F$-block. The following theorem connects exotic and block-exotic fusion systems.

\begin{theorem} (Brauer's Third Main Theorem)
\label{indiana}
\cite[Theorem 7.1]{radha} \textit{Let $G$ be a finite group and $b$ the principal block of $kG$ with maximal $b$-Brauer pair $(P,e_P)$. Then, for any $Q \leq_p G$, $\br^G_Q(b)$ is the principal block of $kC_G(Q)$. In particular, $P \in \syl_p(G)$ and $\mathcal{F}_{(P,e_P)}(G,b)=\mathcal{F}_P(G)$}.
\end{theorem}

In particular, any realisable fusion system is block-realisable. The converse is still an open problem, as noted in Conjecture \ref{XX}.

\subsection{Generalised block fusion systems}
\label{general}

We need to introduce more general categories than block fusion systems, since some group theoretic properties are not captured by these: Assume $b$ is a block of $kG$ with maximal $b$-Brauer pair $(P,e_P)$ and $N \unlhd G$. If $c$ is a block of $kN$ covered by $b$, i.e. $bc \neq 0$, $P \cap N$ is a defect group for $c$, see \cite[Chapter 5, Theorem 5.16 (iii)]{japan}. However, in general $\mathcal F_{(P \cap N, e_{P \cap N})}(N,c)$ is not even a subsystem of $\mathcal F_{(P,e_P)}(G,b)$, see \cite[Example 7.5]{radha}.\\

We use a generalised category, introduced in \cite{ks}, to circumvent this difficulty, which will turn out to be very useful when proving results about block fusion systems.

\begin{defi}
\textit{Let $G$ be a finite group, $N \unlhd G$ and $c$ be a $G$-stable block of $kN$. A \textbf{$(c,G)$-Brauer pair} is a pair $(Q, e_Q)$, where $Q$ is a $p$-subgroup of $G$ with $\br^N_Q(c) \neq 0$ and $e_Q$ is a block of $kC_N(Q)$ such that $\br^N_Q(c)e_Q \neq 0$}.
\end{defi}

Let $(Q,e_Q)$ and $(R,e_R)$ be two $(c,G)$-Brauer pairs. We say that $(Q,e_Q)$ is contained in $(R,e_R)$ and write $(Q,e_Q) \leq (R,e_R)$, if $Q \leq R$ and for any primitive idempotent $i \in (kN)^R$ with $\br^N_R(i)e_R \neq 0$, we also have $\br^N_Q(i)e_Q \neq 0$. This defines an order relation on the set of $(c,G)$-Brauer pairs compatible with the conjugation action of $G$. We also have that given a $(c,G)$-Brauer pair and $Q \leq R$ there exists a unique $(c,G)$-Brauer pair $(Q,e_Q)$ contained in $(R,e_R)$, see \cite[Theorem 1.8(i)]{bp2021}. Also, by \cite[Theorem 1.14(2)]{bp2021}, all maximal $(c,G)$-Brauer pairs are $G$-conjugate. If $(P,e_P)$ is a maximal $(c,G)$-Brauer pair and $G=N$, $P$ is a defect group as defined in the previous section.

\begin{theorem}
\cite[Theorem 3.4]{ks} \textit{Let $G$ be a finite group, $N \unlhd G$ and $c$ be a $G$-stable block of $kN$. Let $(P,e_P)$ be a maximal $(c,G)$-Brauer pair. The category $\mathcal F_{(P,e_P)}(G,N,c)$ is a fusion system. If $(P',e'_P)$ is another maximal $(c,G)$-Brauer pair, then $\mathcal F_{(P,e_P)}(G,N,c)$ is isomorphic to $\mathcal F_{(P',e'_P)}(G,N,c)$}.
\end{theorem}

If $G=N$, generalised block fusion systems coincide with the block fusion systems from Theorem \ref{domingo}. The following theorem shows how this category is connected to the block fusion systems defined before.

\begin{theorem}
\label{liechtenstein}
\cite[Theorem 3.5]{ks} \textit{Let $G$ be a finite group with $N \unlhd G$. Let $c$ be a $G$-stable block of $kN$ covered by a block $b \in kG$. Let $(P,e_P)$ be a maximal $b$-Brauer pair. Then there exists a maximal $(c,G)$-Brauer pair $(S,e_S')$ with $P \leq S$ and $\mathcal F_{(P,e_P)}(G,b) \leq \mathcal F_{(S,e_S')}(G,N,c)$. Further, $P \cap N=S \cap N$, $(S \cap N,e_{S \cap N}')$ is a maximal $c$-Brauer pair and $\mathcal F_{(S \cap N, e_{S \cap N}')}(N,N,c)=\mathcal F_{(S \cap N, e_{S \cap N}')}(N,c) {\unlhd} \mathcal F_{(S,e_S')}(G,N,c)$.}
\end{theorem}

Note that in \cite{ks}, only weak normality of $\mathcal F_{(S \cap N, e_{S \cap N}')}(N,c)$ in $\mathcal F_{(S,e_S')}(G,N,c)$ was proved, but it was improved to normality in \cite[Part IV, Theorem 6.4]{ako}. This theorem shows that the generalised block fusion system contains both the fusion system of the covering block as well as the fusion system of the covered block. Thus these containments form a triangle and we refer to the relations between these fusion systems as ``triangle relations". These useful relations allow for a descent to normal subgroups, which is not possible for block fusion systems as noted earlier. This gives reason to believe that in order to prove Conjecture \ref{XX}, one must prove Conjecture \ref{XZ}.

Note that Conjecture \ref{XZ} implies Conjecture \ref{XX}, since any block fusion system $\mathcal F_{(P,e_P)}(G,b)$ can be written as a generalised block fusion system $\mathcal F_{(P,e_P)}(G,G,b)$.\\

We finish this section by proving that when specialising to the case when $G/N$ is a $p'$-group, then the three fusion systems in Theorem \ref{liechtenstein} are categories on the same group.

\begin{prop}
\label{guernsey}
\textit{Let $G$ be a finite group with normal subgroup $N$ and let $c$ be a block of $kG$ covering a $G$-stable block $b$ of $kN$. If $G/N$ is a $p'$-group, then three fusion systems occurring in Theorem \ref{liechtenstein} are categories on the same group}.
\end{prop}

\textit{Proof}. Let $(S, f_S)$ be a maximal $c$-Brauer pair. As $bc=c$ and $\br^G_S(c)f_S=f_S$, there is a central primitive idempotent $f'_S$ of $kC_N(S)$ such that $\br^N_S(b)f'_S=f'_S$ and $f_S$ covers $f'_S$. In particular, $(S,f'_S)$ is a $(b,G)$-Brauer pair. Let $(T,f'_T)$ be a maximal $(b,G)$-Brauer pair such that $(S,f'_S) \leq (T,f'_T)$. By \cite[Theorem 3.5]{ks}, $(T \cap N, f'_{T \cap N})$ is a maximal $b$-Brauer pair. Since $G/N$ is a $p'$-group we have that $T=T \cap N$ and thus $S \leq T \cap N$. But we also have $|S| \geq |T \cap N|$ which implies $S=T \cap N=T$. \hfill$\square$\\

Note that this proof partially follows the proof of Theorem \ref{liechtenstein}.

\section{Reduction Theorems for generalised block fusion systems}
\label{drei}
Before proving the generalised versions of the Fong Reductions, we set up notation for and state the original Fong Reductions.

\begin{defi}
\textit{Let $H$ be a finite group, $K \leq H$ and $A$ an $H$-algebra over $k$. We denote by $A^H$ the \textbf{$H$-fixed point subalgebra} of $A$.\\
$(a)$ The \textbf{relative trace map} is defined by $\trr^H_K: A^K \rightarrow A^H, a \mapsto {\sum\limits_{x \in [H/K]}} ^x a$, where $[H/K]$ denotes a set of representatives of the right cosets of $K$ in $H$.\\
$(b)$ Denote by $A^H_K$ the image of $A^K$ under the relative trace map $\trr^H_K$.}
\end{defi}

When studying the relations between a group $G$ with a block $b$ and the blocks of a normal subgroup $N \unlhd G$, a certain subgroup of $G$ plays a huge role:

\begin{theorem}
\label{fongone}
\cite[Part IV, Proposition 6.3]{ako} \textit{Let $G$ be a finite group with $N \unlhd G$ and let $c$ be a block of $kN$. Let $I_G(c)=\{g \in G \mid ^gc=c \}$, the stabiliser of $c$ in $G$. The map $e \mapsto \trr^G_{I_G(c)}(e)$ is a bijection between the set of blocks of $kI_G(c)$ covering $c$ and the set of blocks of $kG$ covering $c$. Furthermore, if $b$ is a block of $kG$ covering $c$, then the fusion system of the block $b$ is isomorphic to the fusion system of the block $\tilde b$ of $I_G(c)$ with $b=\trr^G_{I_G(c)}(\tilde b)$}.
\end{theorem}

This theorem is called the First Fong Reduction.

\begin{defi}
\textit{Let $N \unlhd G$, c a block of $kN$ covered by $b \in kG$. We call the block $\tilde b$ of $I_G(c)$ with $b=\trr^G_{I_G(c)}(\tilde b)$ the \textbf{Fong correspondent} of $b$.}
\end{defi}

We use Theorem \ref{fongone} often in the following form:

\begin{cor}
\label{salted}
\textit{Let $\mathcal F$ be a fusion system and $G$ be a finite group possessing an $\mathcal F$-block $b$ such that $|G:Z(G)|$ is minimal among all finite groups having an $\mathcal F$-block. Then $b$ is inertial, i.e. it covers only $G$-stable blocks}.
\end{cor}

\textit{Proof}. Choose $N, c$ as in Theorem \ref{fongone}. Since $Z(G) \subseteq I(c)$, this proposition implies directly that $b$ is inertial. \hfill$\square$\\

The following theorem is the Second Fong Reduction:

\begin{theorem}
\label{fongtu}
\cite[Part IV, Theorem 6.6]{ako} \textit{Let $G$ be a finite group with $N \unlhd G$ and $c$ be a $G$-stable block of $kN$ with trivial defect. Let $b$ be a block of $kG$ covering $c$ and let $(P,e_P)$ be a maximal $b$-Brauer pair. Then $N \cap P=1$ and there exists a central extension $1 \rightarrow Z \rightarrow \widetilde{G} \rightarrow G/N \rightarrow 1$, where $Z$ is a cyclic $p'$-group such that there is a block $\widetilde b$ of $k \widetilde G$ such that if we identify $P$ with the Sylow $p$-subgroup of the inverse image of $PN/N$ in $\widetilde G$, then there is a maximal $\widetilde b$-Brauer pair $(P,f_P)$ such that $\mathcal F_{(P,e_P)}(G,b)=\mathcal F_{(P,f_P)}(\widetilde G, \widetilde b)$}.
\end{theorem}

Now we are to prove new results for generalised block fusion systems. We start with Brauer's Third Main Theorem.\\

\textit{Proof of Theorem \ref{a}}. Let $b$ and $d$ be the principal blocks of $kG$, respectively of $kH$. Then $db \neq 0$ and $b$ is $H$-stable. By Brauer's Third Main Theorem (see Theorem \ref{indiana}), $\br^H_Q(d)$ is the principal block of $kC_H(Q)$ for any $Q \leq_p H$. In particular, there is a maximal $d$-Brauer pair $(P,e_P)$, where $P \in \syl_p(H)$ and $e_P$ is the principal block of $kC_H(P)$ and $\mathcal F_{(P,e_P)}(H,d)=\mathcal F_P(H)$. By Theorem \ref{liechtenstein}, there exists a maximal $(b,H)$-Brauer pair $(S,e_S')$ such that $\mathcal F_P(H)=\mathcal F_{(P,e_P)}(H,d) \subseteq \mathcal F_{(S,e_S')}(H,G,b) \subseteq \mathcal F_S(H)$. Since $P \in \syl_p(H)$, it follows that $S=P$ and that $\mathcal F_{(S,e_S')}(H,G,b)=\mathcal F_P(H)$.\\
By the above, $S=P$. For any $Q \leq P$, let $e_Q$ be the unique block of $kC_H(Q)$ such that $(Q,e_Q) \leq (P,e_P)$ as $d$-Brauer pairs and let $e_Q'$ be the unique block of $kC_G(Q)$ such that $(Q,e_Q') \leq (P,e_P')$ as $(b,H)$-Brauer pairs. By Lemma \ref{rtrt}, $e_Qe_Q' \neq 0$. In other words, the block $e_Q$ of $kC_H(Q)$ covers the block $e_Q'$ of $kC_G(Q)$. But $e_Q$ is the principal block of $kC_H(Q)$. Clearly, the intersection of a Sylow $p$-subgroup and a normal subgroup of any finite group is a Sylow $p$-subgroup of this normal subgroup. This implies that a principal block covers only principal blocks. In particular, $e_Q'$ is the principal block of $kC_G(Q)$. \hfill$\square$\\

Note that if $G=H$, Theorem \ref{a} becomes the original Brauer's Third Main Theorem (see Theorem \ref{indiana}).\\

Next, we generalise both Fong reductions, starting with the first one. We need some background before proving it. We follow the approach of \cite[Chapter 8.7]{block} and directly quote some of the background needed. Fix $(K,\mathcal O,k)$ to be a $p$-modular system.

\begin{defi}
\textit{Let $G$ be a finite group, $A$ be a $G$-algebra over $\mathcal O$ and $P \leq_p G$.\\
$(a)$ We set 
\begin{equation*}
A(P)=A^P/(\sum\limits_{Q<P} A^P_Q + J(\mathcal O)A^P)
\end{equation*} and call the canonical map $\br^A_P: A^P \rightarrow A(P)$ the \textbf{Brauer homomorphism}.\\
$(b)$ Let $H \leq G$. A \textbf{point} of $H$ on $A$ is an $(A^H)^\times$-conjugacy class $\gamma$ of primitive idempotents in $A^H$.\\
$(c)$ A \textbf{local point} of $P$ on $A$ is a point $\gamma$ of $P$ on $A$ with $\br^A_P(\gamma) \neq 0$.}
\end{defi}

Let $G, H$ be finite groups with $G \unlhd H$ and $b$ be an $H$-stable block of $kG$ with $P \leq_p H$ being maximal such that $\br^G_P(b) \neq 0$. A source idempotent is a primitive idempotent $i \in (kGb)^P$ satisfying $\br^G_P(i) \neq 0$ such that for any $Q \leq P$, there is a unique block $e_Q$ with $\br^G_Q(i)e_Q \neq 0$. The interior $P$-algebra $A=i k H i$ is called a source algebra of the block $b$. The pair $(P,e_P)$ is thus a maximal $(b,H)$-Brauer pair. In particular, any source idempotent $i$ determines a generalised block fusion system $\mathcal F=\mathcal F_{(P,e_P)}(H,G,b)$ on $P$ and we call $\mathcal F$ the fusion system of the source algebra $A$ or the fusion system on $P$ determined by the source idempotent $i$.\\
Let $R$ be $\mathcal O$ or  $k$. For an $RG$-module $M$ and a group homomorphism $\varphi: G \rightarrow R^\times$, denote by $_{\varphi} M$ the $RG$-module obtained from  restricting the action of $RG$ along $\varphi$, that is, $_{\varphi} M=M$ as $R$-module and $g \in G$ acts on $m \in M$ as a multiplication by $g\varphi(g)^{-1}$.

\begin{prop}
\label{one}
\textit{Let $G \unlhd H$, $b$ an $H$-stable block of $kG, (P,e_P)$ a maximal $(b,H)$-Brauer pair, $i$ a source idempotent of $(kGb)^P$ such that $\br^G_P(i)e_P=\br^G_P(i)$. Then $\mathcal F_{(P,e_P)}(H,G,b)$ is generated by the set of inclusions between subgroups of $P$ and those automorphisms $\varphi$ of subgroups $Q$ of $P$ such that $_{\varphi} {kQ}$ is isomorphic to a direct summand of $ikHi$ as $k[Q \times Q]$-module}.
\end{prop}

We need two auxiliary results before proving this proposition.

\begin{prop}
\textit{Let $G \unlhd H$, $b$ an $H$-stable block of $kG, (P,e_P)$ a maximal $(b,H)$-Brauer pair, $i$ a source idempotent of $(kGb)^P$ and $A=i kH i$. Let $\mathcal F$ be the fusion system of $A$ on $P$ and $Q \leq P$ fully $\mathcal F$-centralised with $e_Q$ being the unique idempotent such that $(Q,e_Q) \leq (P,e_P)$. For any local point $\delta$ of $Q$ on $\mathcal OHb$ satisfying $\br^G_Q(\delta)e_Q \neq \{0\}$, we have $\delta \cap A \neq \emptyset$}.
\end{prop}

\textit{Proof}. By definition, if $\delta$ is a local point of $Q$ on $\mathcal OHb$ satisfying $\br^G_Q(\delta)e_Q \neq \{0\}$, then $\br^G_Q(\delta)$ is a conjugacy class of primitive idempotents in $kC_H(Q)e_Q$. Since $kC_H(Q)e_Q$ is Morita equivalent to $A(Q)=\br^G_Q(i)kC_H(Q)\br^G_Q(i)$, see \cite[Theorem 6.4.6]{block}, it follows that there is $j \in \delta$ such that $\br^G_Q(j) \in \br^G_Q(i)kC_H(Q)\br^G_Q(i)$. The lifting theorems for idempotents, see \cite[Theorem 4.7.1]{block}, imply that $j$ can be chosen in $A^Q=i(kH)^Qi$. \hfill$\square$

\begin{prop}
\textit{Keep the assumptions of Proposition \ref{one}, let $A=i kH i$, $\mathcal F$ be the fusion system of $A$ on $P$ and $Q,R \leq P$.\\
$(i)$ Every indecomposable direct summand of $A$ as $\mathcal OQ$-$\mathcal OR$-bimodule is isomorphic to $\mathcal OQ \otimes_\varphi \mathcal O R$ for some $ S \leq Q$ and $\varphi \in \homo_{\mathcal F}(S,R)$;\\
$(ii)$ If $\varphi \in \homo_{\mathcal F}(S,R)$, $R$ fully $\mathcal F$-centralised, then $_{\varphi} \mathcal O R$ is isomorphic to a direct summand of $A$ as $\mathcal O Q$-$\mathcal O R$-bimodule}.
\end{prop}

\textit{Proof}. $(i)$ Let $Y$ be an indecomposable direct summand of $A$ as an $\mathcal OQ$-$\mathcal OR$-bimodule. Then $Y$ has a $Q \times R$-stable $\mathcal O$-basis on which $Q$ and $R$ act freely on the left and on the right, respectively. So $Y \cong \mathcal OQ \otimes_{\mathcal OS} {_\varphi \mathcal OR}$ for some $S \leq Q$ and some injective group homomorphism $\varphi: S \rightarrow R$. Let $T=\varphi(S)$. Restricting $Y$ to $S \times T$ shows that $_\varphi \mathcal O T$ is isomorphic to a direct summand of $Y$, hence of $A$, as an $\mathcal OS$-$\mathcal OT$-bimodule. Now $A$ is a direct summand of $\mathcal O H$ as an $\mathcal OP$-$\mathcal OP$-bimodule. In particular, $_\varphi \mathcal O T$ is isomorphic to a direct summand of $\mathcal O H$ as an $\mathcal OS$-$\mathcal OT$-bimodule, hence isomorphic to $ \mathcal OSy^{-1}=y^{-1} \mathcal O T$ for some $y \in H$ with $^yS=T$ and $^ys=\varphi(s)$ for all $s \in S$. Then $\mathcal OS$ is isomorphic to a direct summand of $i \mathcal OHiy=i \mathcal OHy^{-1}iy$ as an $\mathcal OS$-$\mathcal OS$-bimodule. Thus $\br^G_{S}(i \mathcal OH y^{-1}iy) \neq 0$ by \cite[Lemma 5.8.8]{block}. Since $\br^G_{S}(i) \in kC_H(S)e_S$ this also forces $\br^G_{S}(y^{-1}iy)e_S \neq 0$. Conjugating by $y$ yields $\br^G_{T}(i)^y e_S \neq 0$. But this means $^y e_S=e_T$ because $e_T$ is the unique block of $kC_H(T)$ with the property $\br^G_{T}(i)e_T \neq 0$. This shows that $\varphi \in \mor(\mathcal F)$, whence $(i)$.\\
$(ii)$ By definition, there is $x \in H$ such that $\varphi(u)$=$^xu$ for all $u \in Q$ and $^xe_Q=e_R$. Let $\mu$ be a local point of $Q$ on $\mathcal OHb$ such that $\mu \cap A \neq \emptyset$. Set $\nu$=$^x \mu$, i.e. $^x(Q_\delta)=R_\nu$. $\mu \cap A \neq \emptyset$ implies $\br^G_{Q}(\mu)e_Q \neq 0$. Conjugating by $x$ gives $\br^G_{R}(\mu)e_R \neq 0$. Since $R$ is fully $\mathcal F$-centralised, we get from the previous proposition that $\nu \cap A \neq 0$. Let $m \in \mu \cap A$ and $n \in \nu \cap A$. Then $^xm$ and $n$ both belong to $\nu$, so they are conjugate in $(A^R)^\times$. Since $\br^G_{Q}(m) \neq 0$, we get $(m \mathcal OHm)(Q) \neq \{0\}$. This means $m \mathcal OHm$ has a direct summand isomorphic to $\mathcal OQ$ as an $\mathcal OQ$-$\mathcal OQ$-bimodule. This means $m \mathcal OHmx^{-1}=m \mathcal O H xmx^{-1} \cong m \mathcal O H m = mAn$ has a direct summand isomorphic to $\mathcal OQ_{\varphi^{-1}} \cong _\varphi \mathcal O R$ as an $\mathcal OQ$-$\mathcal OR$-bimodule. \hfill$\square$\\

\textit{Proof of Proposition \ref{one}}. By Alperin's Fusion Theorem, see Theorem \ref{alperin}, all morphisms in $\mathcal F_{(P,e_P)}(H,G,b)$ are generated by automorphisms of fully $\mathcal F$-centralised subgroups of $P$. Let $Q$ be such a group. If we set $Q=R$ in the previous proposition, we thus obtain Proposition \ref{one}. \hfill$\square$

\begin{prop}
\label{2023}
\textit{Let $N \unlhd G \unlhd H$ with $N \unlhd H$, $c$ be a block of $kN$ covered by $b \in kG$ and let $\tilde b$ be the Fong correspondent of $b$ in $I_G(b)$. There exists a subgroup $P$ of $I_H(c)$ and $i \in (kI_G(c) \tilde b)^P$ such that $i$ is a source idempotent of the $I_H(c)$-algebra $kI_G(c) \tilde b$ and of the $H$-algebra $kGb$}.
\end{prop}

\textit{Proof}. By the Frattini argument, $H=GI_H(c)$ ($G \unlhd H$ and $G$ transitively permutes all $H$-conjugates of $c$).\\
In particular, $H/I_H(c)=GI_H(c)/I_H(c)=G/(G \cap I_H(c))=G/I_G(c)$. Thus, a system of left coset representatives of $I_G(c)$ in $G$ is also a system of left coset representatives of $I_H(c)$ in $H$. Since $\tilde b$ is $I_H(c)$-stable, we obtain $b=\trr^G_{I_G(c)} \tilde b=\trr^H_{I_H(c)} \tilde b$. Let $x \in H \setminus I_H(c)$, then $\tilde b^x \tilde b=\tilde b c^x c^x \tilde b=0$.\\
We are viewing $kGb$ as $H$-algebra and $kI_G(c) \tilde b$ as $I_H(c)$-algebra. Let $P \leq H$ be maximal such that $\br^G_P(b) \neq 0$ and $j \in (kGb)^P$ primitive idempotent such that $\br^G_P(j) \neq 0$. Let $\mathfrak m$ be the maximal ideal of $(kGb)^P$ not containing $j$ and let $\pi: (kGb)^P \rightarrow (kGb)^P/\mathfrak m$ be the canonical surjection. Since $\br^G_P(j) \neq 0$, we have a factorisation, see \cite[Lemma 14.4]{thevenaz},
\begin{center}
\begin{tikzcd} (kGb)^P \arrow[rd, "\br^G_P"]\arrow[r, "\pi"] & (kGb)^P/\mathfrak m \\& \br^G_P((kGb)^P) \arrow[u, "\pi_0"]\end{tikzcd}
\end{center}
We have: $0 \neq \pi(j)=\pi(jb)=\pi(j) \pi(b)$. Hence, $\pi(b) \neq 0$. Further $\pi(b)=\pi_0 \circ \br^G_P(\trr^H_{I_H(c)} \tilde b)=\pi_0 (\br^G_P(\trr^H_{I_H(c)} \tilde b))$. By the Mackey formula, see \cite[Proposition 2.5.5]{block}, the latter expression is equal to $\pi_0(\br^G_P(\sum\limits_{x \in P \backslash H /I_H(c)} \trr^P_{P \cap ^x I_H(c)} (^x  \tilde b)))$, which means that $\pi(b)=\pi_0(\sum\limits_{x \in P \backslash H /I_H(c)} \br^G_P(\trr^P_{P \cap ^x I_H(c)} (^x \tilde b)))$.\\
However, $\br^G_P(\trr^P_R a)=0$ for any proper subgroup $R$ of $P$ and $a \in A^R$. Thus, the above shows that there is $x \in H$ such that $P \cap ^xI_H(c)=P$ and $0 \neq \pi_0(\br^G_P(^x \tilde b))=\pi(^x \tilde b)$.\\
The equation $0 \neq \pi(^x \tilde b)$ implies that there is a primitive idempotent $i$ of $(kGb)^P$ such that $i$ is $(kGb)^P$-conjugate to $j$ and such that $i^x \tilde b i=i$, i.e. $i \in ^x \tilde b kG ^x \tilde b$. In particular, $\pi(i) \neq 0$, hence $\br^G_P(i) \neq 0$. Setting $i'=^{x^{-1}}i, P'=^{x^{-1}}P$ we get: $i'$ is a primitive idempotent of $(kGb)^{P'}$ with $\br^{I_G(c)}_{P'}(i') \neq 0$ and $i' \in \tilde b kG \tilde b=kI_G(c) \tilde b$.\\
Thus, we have shown: If $j$ is a primitive idempotent of $(kGb)^P$ such that $\br^G_P(j) \neq 0$, then there exists $x \in H$ such that $^{x^{-1}}P \subseteq I_H(c)$ and a $(kGb)^{P}$-conjugate $i$ of $j$ such that $^{x^{-1}}i$ is a primitive idempotent of $kI_G(c) \tilde b$ with $\br^{I_G(c)}_{^{x^{-1}}P}(^{x^{-1}}i) \neq 0$.\\
Replacing $(P,i)$ with $(^{x^{-1}}P,^{x^{-1}}i)$ we obtain: There exists a $p$-subgroup $P$ of $I_H(c)$ and a primitive idempotent $i$ of $kI_H(c) \tilde b$ such that $\br^G_P(i) \neq 0$, $i$ is primitive in $kGb$ and $P$ is maximal among $p$-subgroups of $H$ such that $\br^{I_G(c)}_P(b) \neq 0$.\\
Conversely, suppose that $P \leq I_H(c)$ is maximal such that $\br^{I_G(c)}_P(\tilde b) \neq 0$. By maximality, $\tilde b=\trr^{I_H(c)}_P(a)$ for some $a \in (kG \tilde b)^P$. Thus, $b=\trr^H_{I_H(c)} \tilde b=\trr^H_P(a)$. This shows that $P$ is contained in a maximal $p$-subgroup $Q$ of $H$ such that $\br^G_Q(b) \neq 0$.\\
Combining, we obtain: There exists a $p$-subgroup $P$ of $I_H(c)$ and a primitive idempotent $i$ of $(kI_G(c) \tilde b)^P$ such that $P$ is maximal among subgroups $Q$ of $H$ such that $\br^G_Q(b) \neq 0$, $P$ is maximal among subgroups $Q$ of $I_H(c)$ such that $\br^{I_G(c)}_Q(\tilde b) \neq 0$, $i$ is a primitive idempotent of $(kGb)^P$ and $\br^G_P(i) \neq 0$. \hfill$\square$

\begin{prop}
\label{angoracat}
\textit{Keep the notation from the previous proposition. Then the algebras $i  k Hi$ and $i k I_H(c)i$ are isomorphic as interior $P$-algebras}.
\end{prop}

\textit{Proof}. Let $P, i$ be as above, then we claim that also have $ikHi=ikI_H(c)i$. Indeed, clearly $ikI_H(c)i \subseteq ikHi$.\\
Now suppose $x \in H \setminus I_H(c)$. Then $ixi=i \tilde b x \tilde b i=i \tilde b ^x \tilde b xi=0$. This shows that $ikHi \subseteq ikI_H(c)i$. \hfill$\square$\\

\textit{Proof of Theorem \ref{b}}. Since $I_G(c)=I_H(c) \cap G$, and $G$ is normal in $H$, $I_G(c)$ is normal in $I_H(c)$.\\
Next, we claim that $\tilde b$ is $I_H(c)$-stable. Let $x \in I_H(c)$, then
\begin{equation*}
b=^xb=^x \trr^G_{I_G(c)}(\tilde b)=\trr^{^xG}_{^xI_G(c)}(^x \tilde b)=\trr^G_{I_G(c)}(^x \tilde b),
\end{equation*} where the last equality follows from the normality of $G$ in $H$ and $I_G(c)$ in $I_H(c)$. This equation, together with the uniqueness of the Fong correspondent, implies stability. Note that $^x \tilde b$ is indeed a block of $kI_G(c)$ covering $c$ since $^x \tilde b c$=$^x \tilde b ^x c$=$^x(\tilde b c) \neq 0$. \\
We apply Propositions \ref{2023} and \ref{angoracat} to obtain a source idempotent $i$ of $b$ and $\tilde b$ respectively with $i kHi=i  k I_H(c) i$. This observation, together with Proposition \ref{one}, implies the theorem.  \hfill$\square$\\

Theorem \ref{b} is a generalisation of the First Fong Reduction, which we obtain from it if $G$ and $H$ coincide. We give an example proving that the assumption $N \unlhd H$ in the theorem is necessary. 

\begin{bsp}
Let $b$ and $c$ be principal blocks, then, by Theorem \ref{a}, the statement becomes $\mathcal F_S(H)=\mathcal F_S(I_H(c))$. In particular, since $c$ is principal, we also have $I_H(c)=N_H(N)$, thus $\mathcal F_S(H)=\mathcal F_S(N_H(N))$ whenever $N \unlhd G \unlhd H$.\\
Now let $p=3$, $N=(C_3 \times C_3) \rtimes C_2$, where $C_2$ acts as a reflection, then $S=C_3 \times C_3 \in \syl_p(N)$. Let $H=S \rtimes D_8$, $G=S \rtimes (C_2 \times C_2)$ then $N_H(N)=G$, but clearly $\mathcal F_S(G) \subsetneq \mathcal F_S(H)$.
\end{bsp}

Finally, we also generalise the Second Fong Reduction to generalised block fusion systems.\\

\textit{Proof of Theorem \ref{c}}. First, we apply the second Fong reduction to $A \unlhd H$. We use the notation from \cite[Proof of Theorem 3.1]{solom} and are recalling some key steps. Let $S$ be a $p$-subgroup of $H$ containing $P$, maximal such that $\br^A_S(c) \neq 0$. For each $h \in H$ there is an element $i_h \in (kAc)^\times$ such that $c_{i_h}=c_h$ on $(kAc)^\times$. We can choose the elements $i_h$ such that $i_{ha}=i_{h}ac$ for $h \in H, a \in A$ and such that
\begin{equation*}
S \rightarrow (kAc)^\times, s \mapsto i_s
\end{equation*}
is a homomorphism. Now define a 2-cocycle $\alpha$ on $H$ via $i_gi_h=\alpha(g,h)i_{gh}$ for $g,h \in H$. Denote by $k_{\overline{\alpha}^{-1}}{H/A}$ the twisted group algebra corresponding to $\overline \alpha^{-1}$, i.e. the free module on $\{\widehat{\overline h} \mid \overline h \in H/A\}$ with multiplication given by $\widehat{\overline h} \widehat{\overline g}=\alpha^{-1}(h,g)\widehat{\overline{hg}}$. Define a function $\phi: kAc \otimes k_{\overline{\alpha}^{-1}}\widehat{H/A} \rightarrow kHc$ via the $k$-linear extension of the map $x \otimes \widehat{\overline h} \mapsto xi_{h}^{-1}h$. This gives a central $p'$-extension
\begin{equation*}
1 \rightarrow Z \rightarrow \widetilde H \rightarrow H/A \rightarrow 1
\end{equation*}
and by our choice of the elements $i_h$ also a $p'$-extension
\begin{equation*}
1 \rightarrow Z \rightarrow \widetilde M \rightarrow M/A \rightarrow 1,
\end{equation*} where $\widetilde M$ is the full inverse image of $M/A$ in $\widetilde H$. Furthermore, we get an idempotent $e$ of $kZ$ and an algebra isomorphism $\tau: k \widetilde H e \rightarrow k_{\overline{\alpha}^{-1}}{\widehat{H/A}}$ with $\tau(\tilde he)=\alpha_{\tilde h} \widehat{\overline h}$ for some $\alpha_{\tilde h} \in k^\times$, see \cite[Proposition 10.8]{thevenaz}. Denote by $\tau_M$ the restriction to $k \widetilde M e$. If we define $\alpha_M$ to be the restriction of $\alpha$ to $M \times M$, then $\tau_M$ becomes an algebra isomorphism from $k \widetilde M e$ to $k_{\overline{\alpha_M}^{-1}}\widehat{M/A}$.\\
For each $s \in S$, let $\tilde s$ denote the unique lift of $s$ in $\widetilde H$ which is also a $p$-element and for $Q \leq S$ define $\widetilde Q=\{\tilde q \mid q \in Q\}$. Note that the groups $S/A, \widetilde S$ and $S$ are isomorphic and we identify them henceforth. In particular, $\widetilde P \cong P$. If we consider $kAc \otimes k\widetilde H e$ as interior $S$-algebras via $s \mapsto i_s \otimes se$, we obtain an $S$-algebra isomorphism
\begin{equation*}
 \psi: kAc \otimes k\widetilde H e \rightarrow kHc, x \otimes y \mapsto \phi(x \otimes \tau(y)).
\end{equation*}
Again, we denote the restriction to $kAc \otimes k\widetilde M e$, which is also an $S$-algebra isomorphism to $kMc$, by $\psi_M$. Since $kAc$ and $k \widetilde H e$ are $p$-permutation algebras, see \cite[Proposition 28.3]{thevenaz}, $\psi$ induces algebra isomorphisms $\psi_Q: kAc(Q) \otimes k\widetilde H e(Q) \rightarrow kHc(Q)$ and $\psi_{Q,M}: kAc(Q) \otimes k\widetilde M e(Q) \rightarrow kMc(Q)$ satisfying $\psi_{Q,M}(\br^M_Q(x) \otimes \br^M_Q(y))=\br^M_Q(\psi_M(x \otimes y))$. Since $kAc$ is a matrix algebra, we get bijections between the blocks of $H$ covering $c$ and the blocks of $\widetilde H$ covering $e$ as well as between the blocks of $M$ covering $c$ and the blocks of $\widetilde M$ covering $e$. For the same reason, $\psi_{Q,M}$ induces a bijection between the blocks of $kMc(Q)$ and $k\widetilde Me(Q)$ via $f \mapsto \psi_{Q,M}(1 \otimes f)$. Thus, $(Q,f) \mapsto (Q,\psi_{Q,M}(1 \otimes f))$ provides a bijection between the set of Brauer pairs associated to blocks of $\widetilde M$ covering $e$ and the set of Brauer pairs associated to blocks of $M$ covering $c$ with first components respectively contained in $S$. \\
Let $\tilde d$ be the block of $k \widetilde M$ corresponding to $d$ under $\psi_M$, i.e. $d=\psi_M(c \otimes \tilde d)$. Note that $d$ is $H$-stable if and only if $d \in Z(kHc)$. Since $\psi_M$ is an $k$-algebra isomorphism, this is if and only if $c \otimes \tilde d \in Z(k \widetilde H \tilde c \otimes kAc)$. The latter expression is equal to $Z(k \widetilde H \tilde c) \otimes k$. Since $d=\psi_M(c \otimes \tilde d)$, we get that $d$ is $H$-stable if and only if $\tilde d$ is $\widetilde H$-stable.\\
Define $\widetilde{e_P}$ by $\psi_{M,P}(c \otimes {\widetilde{e_P}})=e_P$. Note that by the description above, $\psi_{M,P}$ is an inclusion-preserving map of Brauer pairs. Since $P \leq M$ we thus get that $(P,\widetilde{e_P})$ is a maximal $(\tilde d, \widetilde H)$-Brauer pair. We can thus define the category $\mathcal F_{(P,\widetilde{e_P})}(\widetilde H,\widetilde M,\tilde d)$. Now we define a map from $\mathcal F_{(P,e_P)}(H,M,d)$ to $\mathcal F_{(P,\widetilde{e_P})}(\widetilde H,\widetilde M,\tilde d)$ by $(Q,e_Q) \mapsto (Q, \psi_{M,Q}(1 \otimes e_Q))$. The objects of $\mathcal F_{(P,e_P)}(H,M,d)$ and $\mathcal F_{(P,\widetilde{e_P})}(\widetilde H,\widetilde M,\tilde d)$ are the same by the above. For the morphisms we want to prove that $N_H(Q,e_Q)$ and $N_{\widetilde H}(Q,\psi_{M,Q}(1 \otimes e_Q))$ have the same image in $\aut(Q)$. For this, we need to prove that $\psi_{M,Q}$ is $H$-equivariant. Note that we can apply \cite[3.4]{solom} to get that for any $h \in H$ with lift $\tilde h$ of $hA$ in $N_{\widetilde H}(Q)$, there exists some $a \in A$ such that $n:=ah \in N_H(Q)$ and $\widetilde{^nx}=^{\tilde h}\tilde x$ for all $x \in Q$. Let $\widetilde{e_Q}:=\psi_{M,Q}(1 \otimes e_Q)$ and $\widetilde{f_Q} \in (k\widetilde Me)^Q$ such that $\br_Q^M(\widetilde{f_Q})=\widetilde{e_Q}$. If we define $f_Q:=\psi_M(1 \otimes \widetilde{f_Q})$, then $\br_Q^M(f_Q)=e_Q$ by the observations about $\psi_M$ above. Let $h \in N_H(Q)$ with $\tilde h \in \widetilde H$ a lift of $\overline h$. If we define $t:=\psi^{-1}(hc)$, then $t=i_h \otimes \alpha \tilde h e$ for some $\alpha \in k^\times$. We get:
\begin{align*}
& ^he_Q=^h\br^M_Q(f_Q)=\br^M_Q(^hf_Q)=\br^M_Q(\psi(^t(1 \otimes \widetilde{f_Q}))) =
\br^M_Q(\psi(1 \otimes ^{\tilde h}\widetilde{f_Q})) \\
&= \br^M_Q(\psi_M(1 \otimes ^{\tilde h}\widetilde{f_Q})) =\psi_{M,Q}(1 \otimes \br^M_Q(^{\tilde h}\widetilde{f_Q}))=\psi_{M,Q}(1 \otimes ^{\tilde h}\br^M_Q(\widetilde{f_Q})) \\
& =\psi_{M,Q}(1 \otimes ^{\tilde h}\widetilde{e_Q})
\end{align*}
In particular, $h$ stabilises $e_Q$ if and only if $\tilde h$ stabilises $\widetilde{e_Q}$. Thus, $N_H(Q,e_Q)$ and $N_{\widetilde H}(Q,\psi_{M,Q}(1 \otimes e_Q))$ have the same image in $\aut(Q)$ as claimed. By Theorem \ref{alperin}, the morphisms in the categories $\mathcal F_{(P,e_P)}(H,M,d)$ and $\mathcal F_{(P,\widetilde{e_P})}(\widetilde H,\widetilde M,\tilde d)$ coincide, which implies the theorem. \hfill$\square$\\

If $H=M$, we obtain the original Second Fong Reduction.

\section{Extension of Cabanes' Theorem}
\label{funf}
In this section we prove Theorem \ref{kombucha}. We start by giving some background on  character and block theory of such groups. See \cite{oh} for details.

\subsection{Background}
Let $\mathbf G$ be a connected, reductive, algebraic group defined over $\overline{\mathbb{F}}_q$ where $q$ is a prime power and let $F: \mathbf G \rightarrow \mathbf G$ be a Frobenius morphism defining an $\mathbb{F}_q$-structure on $\mathbf G$. Let $\mathbf L \leq \mathbf G$ be an $F$-stable Levi subgroup of a parabolic subgroup $\mathbf P$ of $\mathbf G$.\\
Recall the Lusztig induction map 
\begin{equation*}
R^{\mathbf G}_{\mathbf L}: \mathbb Z \irr(\mathbf L^F) \rightarrow \mathbb Z \irr(\mathbf G^F), [M] \mapsto \sum\limits_{i \in \mathbb Z} (-1)^i [H^i_c(\mathbf Y_{\mathbf P}) \otimes_{\mathbb C \mathbf L^F}M],
\end{equation*}
where $[M]$ is the class of a $\mathbb C \mathbf L^F$-module $M$ in $\mathbb \irr(\mathbf L^F)$, $\mathbf Y_{\mathbf P}$ is the variety 
\begin{equation*}
\{gR_u(\mathbf P) \mid g^{-1}F(g) \in R_u(\mathbf P)F(R_u(\mathbf P)) \},
\end{equation*}
 and $R_u(\mathbf H)$ is the unipotent radical of an algebraic group $\bf H$. This is defined using a parabolic subgroup $\mathbf P$ of which $\mathbf L$ is a Levi. However, in all cases we consider, $R^{\mathbf G}_{\mathbf L}$ is independent of the choice of $\mathbf P$ and hence we will suppress it from the notation. We use this construction in the case where $\mathbf L$ is a maximal torus to parametrise the irreducible characters of $\mathbf G^F$.\\

If we fix an $F$-stable maximal torus of $\mathbf G$, we can define $\mathbf G^\ast$ to be a group in duality with $\mathbf G$ with respect to this torus, with corresponding Frobenius again denoted by $F$. Then there is a bijection $\{\mathbf G^F$-conjugacy classes of pairs $(\mathbf T,\theta)$ where $\mathbf T$ is an $F$-stable maximal torus of $\mathbf G$ and $\theta \in \irr(\mathbf T^F)\} \leftrightarrow \{\mathbf G^{\ast^{F}}$-classes of pairs $(\mathbf T^\ast,s)$ where $\mathbf T^\ast$ is an $F$-stable maximal torus of $\mathbf G^\ast$ and $s \in \mathbf T^{\ast^{F}}\}$.

\begin{theorem}
(Deligne--Lusztig). \cite[Theorem 4.7]{oh} \textit{For $s \in \mathbf G^{\ast^{F}}$ a semisimple element, one defines $\mathcal E(\mathbf G^F,s)$ to be the set of irreducible components of generalised characters $R^{\mathbf G}_{\mathbf T} \theta$ for $(\mathbf T,\theta)$ corresponding to some $(\mathbf T^\ast,s)$ through the above correspondence. One gets a partition $\irr (\mathbf G^F)=\bigsqcup\limits_{s} \mathcal E(\mathbf G^F,s)$, where $s$ ranges over representatives of semisimple classes of $\mathbf G^{\ast^{F}}$. This partition is called Jordan decomposition of $\irr(\mathbf G^F)$}.
\end{theorem}

This theorem plays an important role when studying blocks of groups of Lie type. The blocks obtained when specialising to $s=1$ are the starting point for many constructions.

\begin{defi}
\textit{Keep the notation of the previous theorem. The elements of $\mathcal E(\mathbf G^F,1)$ are called \textbf{unipotent characters}. A block $b$ of $\mathbf G^F$ such that $\irr(\mathbf G^F,b) \cap \mathcal E(\mathbf G^F,1) \neq \emptyset$ is called a \textbf{unipotent block} of $\mathbf G^F$}.
\end{defi}

If $\chi$ is a unipotent character of $\mathbf L^F$, then all constituents of $R^{\mathbf G}_{\mathbf L}\chi$ are unipotent. We recall an important result related to the parametrisation by Deligne--Lusztig.

\begin{theorem} 
\label{naomiwatts}
(Broué--Michel, Hiss). \cite[Theorem 9.12]{caen}
\textit{Assume $(p,q)=1$. Let $s \in \mathbf G^\ast$ be a semisimple $p'$-element, define $b{(\mathbf G^F,s)}=\sum\limits_{\chi \in \mathcal E(\mathbf G^F,s)} e_{\chi}$ and $\mathcal E_p(\mathbf G^F,s)$ as the union of rational series $\mathcal E(\mathbf G^F,t)$ such that $s=t_{p'}$.\\
$(a)$ The set $\mathcal E_p(\mathbf G^F,s)$ is a union of blocks $\irr(\mathbf G^F, b_i)$, i.e. $b(\mathbf G^F,s) \in \mathcal O \mathbf G^F$. \\
$(b)$ For each block $b$ of $k \mathbf G^F$ with $\irr(\mathbf G^F,b) \subseteq \mathcal E_p(\mathbf G^F,s)$, one has $\irr(\mathbf G^F,b) \cap \mathcal E(\mathbf G^F,s) \neq \emptyset$.}
\end{theorem}

If $e \geq 1$, recall that $\phi_e(x) \in \mathbb Z[x]$ denotes the $e$-th cyclotomic polynomial, whose complex roots are the roots of unity of order $e$. Any $F$-stable torus $\mathbf S$ of $\mathbf G$ has so-called polynomial order $P_{\mathbf S,F} \in \mathbb Z[x]$ defined by $|\mathbf S^{F^m}|=P_{\mathbf S,F}(q^m)$ for some $a \geq 1$ and any $m \in 1+a \mathbb N$. Moreover, $P_{\mathbf S,F}$ is a product of cyclotomic polynomials $P_{\mathbf S,F}=\prod\limits_{e \geq 1} \phi_e^{n_e}, n_e \geq 0$. A $\phi_e$-torus of $\mathbf G$ is an $F$-stable torus whose polynomial order is a power of $\phi_e$. An $e$-split Levi subgroup is any $C_\mathbf G(\mathbf S)$, where  $\mathbf S$ is a $\phi_e$-torus of $\mathbf G$.\\
Let $\mathbf L_i, i=1,2,$ be an $e$-split Levi subgroup in $\mathbf G$ and $\zeta_i \in \mathcal E(\mathbf L_i^F,1)$. We write $(\mathbf L_1,\zeta_1) \leq_e (\mathbf L_2, \zeta_2)$ if $\zeta_2$ is a component of $R^{\mathbf L_2}_{\mathbf L_1}(\zeta_1)$, i.e. $\langle R^{\mathbf L_2}_{\mathbf L_1}(\zeta_1), \zeta_2 \rangle \neq 0$.\\
A unipotent character $\chi \in \mathcal E(\mathbf G^F,1)$ is said to be $e$-cuspidal if a relation $(\mathbf G,\chi) \geq_e (\mathbf L, \zeta)$ is only possible with $\mathbf L=\mathbf G$. A pair $(\mathbf L,\zeta)$ with $\mathbf L$ an $e$-split Levi subgroup and an $e$-cuspidal $\zeta \in \mathcal E(\mathbf L^F,1)$ is called a unipotent $e$-cuspidal pair of $\mathbf G^F$.

\begin{theorem} (Cabanes--Enguehard).
\cite[Theorem 6.10]{oh}
\textit{Assume $p \geq 7$ is a prime not dividing $q$. Let $e$ be the multiplicative order of $q$ mod $p$. If $(\mathbf L,\zeta)$ is a unipotent $e$-cuspidal pair of $\mathbf G$, then there exists a block $b_{ \mathbf G^F}(\mathbf L,\zeta)$ such that all constituents of $R^{\mathbf G}_{\mathbf L}(\zeta)$ belong to $b_{ \mathbf G^F}(\mathbf L,\zeta)$. The map $(\mathbf L, \zeta) \mapsto b_{ \mathbf G^F}(\mathbf L,\zeta)$ gives a bijection between $\mathbf G^F$-classes of unipotent $e$-cuspidal pairs of $\mathbf G^F$ and unipotent blocks of $\mathbf G^F$. }
\end{theorem}

A starting point to prove Conjecture \ref{XX} for groups of Lie type in non-describing characteristic is the following theorem by Cabanes, which proves the conjecture for almost all primes if the underlying block is unipotent:

\begin{theorem} (Cabanes) \cite[Theorem 7.11]{oh}
\label{virginia}
\textit{Let $p \geq 7$ a prime with $(p, q) = 1$ and $b$ a unipotent block of
$k\mathbf G^F$ with maximal $b$-Brauer pair $(P,e_P)$. Then the fusion system $\mathcal F_{(P,e_P)}(\mathbf G^F,b)$ is non-exotic.}
\end{theorem}

\subsection{Extension of the Theorem}
We prove Theorem \ref{kombucha} in this subsection. Fix $\mathbf X$ to be a connected, reductive group with Frobenius $F$ defining an $\mathbb{F}_q$-structure on $\mathbf X$, $\mathbf G \unlhd \mathbf N$, where $\mathbf G \leq \mathbf X$ is a Levi subgroup and $\mathbf N \leq N_{\mathbf X}(\mathbf G)$. Let $c$ be a block of $k \mathbf N^F$ covering a unipotent block $b$ of $k\mathbf G^F$ with $b=b_{\mathbf G^F}(\mathbf L,\zeta)$ for an $e$-cuspidal pair $(\mathbf L,\zeta)$ of $\mathbf G^F$.\\

Let $(Q,b_Q)$ be a Brauer pair of a finite group $H$. Then it is called centric if and only if $b_Q$ has defect group $Z(Q)$ in $C_H(Q)$. Note that this is equivalent to $Q$ being centric in the fusion system of a block $b$ of $kH$ such that $(1,b) \leq (Q,b_Q)$. Then there is a single $\zeta \in \irr(b_Q)$ with $Z(Q)$ in its kernel, we call it the canonical character of the centric subpair.\\

We assume for the rest of this section that $p \geq 7$. Let $Z=Z^\circ(\mathbf L)^F_p$. By \cite[Proposition 2.2]{ce93}, $Z=Z(\mathbf L)^F_p$, $\mathbf L=C^\circ_{\mathbf G}(Z)$ and $\mathbf L^F=C^\circ_{\mathbf G}(Z)^F=C_{\mathbf G^F}(Z)$.\\
Let $e_Z$ be the block of $\mathbf L^F$ containing the character $\zeta$. Since $\mathbf L^F=C_{\mathbf G^F}(Z)$, $(Z,e_Z)$ is a Brauer pair for $\mathbf G^F$. By \cite[Lemma 4.5]{ce93}, $(Z,e_Z)$ is a $b$-Brauer pair, $(Z,e_Z)$ is a centric Brauer pair and $\zeta$ is the canonical character of $e_Z$. Further, there exists a maximal $b$-Brauer pair $(P,e_P)$ such that $(Z,e_Z) \unlhd (P,e_P)$ and such that the canonical character of $e_P$ is $\res^{\mathbf L^F}_{C_{\mathbf G^F}(P)}(\zeta)$. Here we note that if $(Q,e)$ is a centric Brauer pair, then for any Brauer pair $(R,f)$ such that $(Q,e) \leq (R,f)$, we also have that $(R,f)$ is also centric.\\
Furthermore, by \cite[Theorem 4.4(ii)]{ce94}, $P$   is a Sylow $p$-subgroup of $(C_{\mathbf G}^{\circ} ([\mathbf L,\mathbf L] ) )^F$.\\

Recall that for a semisimple, algebraic group $\mathbf X$, there exist natural isogenies $\mathbf X_{\simk} \rightarrow \mathbf X \rightarrow \mathbf X_{\ad}$, where $\mathbf X_{\simk}$ is simply connected and $\mathbf X_{\ad}$ is of adjoint type.
There is a decomposition of $[\mathbf X,\mathbf X]=\mathbf X_1 \dots \mathbf X_m$ as a central product of $F$-stable, closed subgroups, see \cite[Definition 22.4]{caen} for more details. Define $\mathbf X_\mathbf a=Z^\circ(\mathbf X)\mathbf X'_\mathbf a$, where $\mathbf X'_\mathbf a$ is the subgroup generated by the $\mathbf X_i$ with $(\mathbf X_i)^F_{\ad} \cong \pgl_{n_i}(\epsilon_iq^{m_i})$ and $p$ dividing $q^{m_i}-\epsilon_i$. Let $\mathbf X_\mathbf b$ be generated by the remaining $\mathbf X_i$. We thus have a decomposition $\mathbf X=\mathbf X_\mathbf a \mathbf X_\mathbf b$.\\

Let $M=N_{\mathbf G^F}([\mathbf L, \mathbf L], \res^{\mathbf L^F}_{[\mathbf L, \mathbf L]^F} \zeta)$. As noted above, $(Z,e_Z) \unlhd (P,e_P)$. Thus, $P$ normalises $\mathbf L=C_{\mathbf G^\circ}(Z)$ and consequently $P$ normalises $[\mathbf L,\mathbf L]$ and $[\mathbf L,\mathbf L]^F$. Again, since $(Z,e_Z) \unlhd (P,e_P)$, $^xe_Z=e_Z$ for all $x \in P$. Since $\zeta$ is the canonical character of $e_Z$, $^x \zeta=\zeta$ for all $x \in P$. Since restriction induces a bijection between $\mathcal E(\mathbf L^F,1)$ and $\mathcal E([\mathbf L,\mathbf L]^F,1)$, see \cite[Proposition 13.20]{dm}, $^x \res^{\mathbf L^F}_{[\mathbf L,\mathbf L]^F} \zeta=\res^{\mathbf L^F}_{[\mathbf L,\mathbf L]^F} \zeta$ for all $x \in P$. Thus, $P \leq M$. Also, $C_{\mathbf G^F}(P) \leq C_{\mathbf G^F}(Z)=\mathbf L^F \leq M$. Since $PC_{\mathbf G^F}(P) \leq M$, $(P,e_P)$ is a Brauer pair for a block of $M$.\\

We recall some details which can be found in the proof of \cite[Theorem 7.11]{oh}: Let $Q \leq P$ be $\mathcal F_{(P,e_P)}(\mathbf G^F,b)$-centric. We have $C^\circ_{\mathbf G}(Q)_{\mathbf b}=[\mathbf L,\mathbf L]$. Let $\zeta_Q$ be the canonical character of $e_Q$. Then $\zeta_Q^\circ:=\res^{C_{\mathbf G^F}(Q)}_{C_{\mathbf G}^\circ(Q)^F} \zeta_Q$ is the unique unipotent character of $C^\circ_{\mathbf G}(Q)^F$ whose restriction to $[\mathbf L,\mathbf L]^F$ is $\res^{\mathbf L^F}_{[\mathbf L,\mathbf L]^F} \zeta$.

\begin{prop}
\label{potas}
\textit{Let $H:=N_{\mathbf N^F}([\mathbf L,\mathbf L], \res^{\mathbf L^F}_{[\mathbf L,\mathbf L]^F} \zeta, d)$, where $d$ is the unique block of $kM$, where $M:=N_{\mathbf G^F}([\mathbf L,\mathbf L], \res^{\mathbf L^F}_{[\mathbf L,\mathbf L]^F} \zeta)$, with $(P,e_P)$ a $d$-Brauer pair. We have
$\mathcal F_{(P,e_P)}(\mathbf N^F,\mathbf G^F,b)=\mathcal F_{(P,e_P)}(H,M,d)$.}
\end{prop}

\textit{Proof}. Let $\mathcal F=\mathcal F_{(P,e_P)}(\mathbf N^F,\mathbf G^F,b)$, $\mathcal G=\mathcal F_{(P,e_P)}(H,M,d)$ and $\mathcal F_0=\mathcal F_{(P,e_P)}(\mathbf G^F,b)$.\\ Assume $(X,e_X) \leq (P,e_P)$ is $\mathcal F_0$-centric. Let $x \in N_{\mathbf N^F}(X,e_X)$, then $x$ normalises $C_{\mathbf G}(X)$, hence, also $C_{\mathbf G}(X)_{\mathbf b}=[\mathbf L,\mathbf L]$. Also, $x$ stabilises $\zeta_X$ and thus also ${\zeta_X}|_{[\mathbf L,\mathbf L]^F}=\res^{\mathbf L^F}_{[\mathbf L,\mathbf L]^F} \zeta$. We thus have $x \in N_{\mathbf N^F}([\mathbf L,\mathbf L]^F,\res^{\mathbf L^F}_{[\mathbf L,\mathbf L]^F} \zeta)$.\\
Furthermore, $C_{\mathbf G^F}(X) \leq M$, so $e_X$ is a block of $kC_M(X)$ and hence $(X,e_X)$ is the unique $d$-Brauer pair with first component $X$ such that $(X,e_X) \leq (P,e_P)$. In particular, $b$-Brauer pairs and $d$-Brauer pairs, for which the first component is a $\mathcal F_0$-centric subgroup of $P$, are the same.\\
Thus, if $x$ stabilises $(X,e_X)$, it also stabilises $d$, and we get $N_{\mathbf N^F}(X,e_X) \leq H$ for all $X \leq P$ which are $\mathcal F_0$-centric. Now if $X \leq P$ is $\mathcal F$-centric, it is in particular $\mathcal F_0$-centric, so $\mathcal F \subseteq \mathcal G$.\\
Conversely, suppose that $(X,e_X) \leq (P,e_P)$ is $\mathcal G$-centric. In particular, by the previous paragraph, $X$ is $\mathcal F$- and thus also $\mathcal F_0$-centric. But again, since $e_X$ is a block of $kC_M(X)$, $(X,e_X)$ is an $(d,H)$-Brauer pair if and only if it is an $(b,\mathbf N^F)$-Brauer pair and hence $\aut_{\mathcal G}(X) \leq \aut_{\mathcal F}(X)$. By Theorem \ref{alperin}, this shows $\mathcal G \subseteq \mathcal F$ and hence $\mathcal G = \mathcal F$. \hfill$\square$\\

Fix $H$ as in the previous proposition for the rest of this section.

\begin{lemma}
\label{mirschmelli}
\textit{For $Q \leq_p \mathbf G^F$, the group $C_{\mathbf G^F}(Q)/C_{\mathbf G}^\circ(Q)^F$ is also a $p$-group.}
\end{lemma}

\textit{Proof}. By \cite[Proposition 2.1.6(e)]{APrairieHomeCompanion}, $C_{\mathbf G}(Q)/C_{\mathbf G}^{\circ} (Q)$ is a $p$-group. Now compose the inclusion $C_{\mathbf G^F}(Q) \rightarrow C_{\mathbf G}(Q)$ with the natural surjection $C_{\mathbf G^F}(Q) \rightarrow C_{\mathbf G^F}(Q)/C_{\mathbf G}^\circ(Q)^F$ and denote this morphism by $\varphi$. Then $\ker \varphi=C_{\mathbf G}^\circ(Q)^F$. Thus $C_{\mathbf G^F}(Q)/C_{\mathbf G}^\circ(Q)^F$ is isomorphic to a subgroup of $C_{\mathbf G}(Q)/C_{\mathbf G}^\circ(Q)$ and in particular also a $p$-group. \hfill$\square$

\begin{lemma}
\label{2020}
\textit{We have $N_{\mathbf N^F}(Q,e_Q)=N_H(Q)$ for any $Q \leq P$ which is $\mathcal F_{(P,e_P)}(\mathbf G^F,b)$-centric.}
\end{lemma}

\textit{Proof}. In the proof of Proposition \ref{potas}, we proved $N_{\mathbf N^F}(Q,e_Q) \leq H$. Thus, in particular $N_{\mathbf N^F}(Q,e_Q) \subseteq N_H(Q)$. Let $x \in N_H(Q)$. We need to prove that $x$ stabilises $e_Q$. First, we prove that $x$ preserves the set of unipotent characters of $C_{\mathbf G}^\circ(Q)^F$. Note that $\mathbf N \subseteq N_{\mathbf X}(\mathbf G)$. Further, $N_{\mathbf N}(Q) \subseteq N_{\mathbf X}(\mathbf G) \cap N_{\mathbf X}(Q) \subseteq N_{\mathbf X}(\mathbf G) \cap N_{\mathbf X}(C_{\mathbf X}(Q)) \subseteq N_{\mathbf X}(C_{\mathbf G}(Q)) \subseteq N_{\mathbf X}(C_{\mathbf G}^\circ(Q))$, whereas the last inclusion follows from the fact that $\mathbf X$ is an algebraic group. Now $x$ normalises $[\mathbf L,\mathbf L], \res^{\mathbf L^F}_{[\mathbf L,\mathbf L]^F} \zeta$ and $C_{\mathbf G}^\circ(Q)$ and thus, by the remarks before Proposition \ref{potas} together with the observation about unipotent characters, also $\zeta^\circ_Q$. In particular, it also stabilises the block of $kC^\circ_{\mathbf G}(Q)^F$ containing $\zeta^\circ_Q$. Further, by Lemmas \ref{mirschmelli} and \ref{feit}, $e_Q$ is the unique block of $kC_{\mathbf G^F}(Q)$ covering the block of $kC_{\mathbf G}^\circ(Q)^F$ containing $\zeta^\circ_Q$. Thus, $x$ stabilises $e_Q$, so $N_H(Q) \subseteq N_{\mathbf N^F}(Q,e_Q)$. \hfill$\square$

\begin{lemma}   
\label{2000}
\textit{Let $A=[\mathbf L,\mathbf L]^F$. The group $PA/ A  $ is in $\syl_p(H/A)$}.
\end{lemma}

\textit{Proof}.  Since $|H:M|_{p} =1 $, it suffices to prove that  $ PA/A$ is a Sylow $p$-subgroup  of $M/A $.    Now  $P$   is a Sylow $p$-subgroup of $(C_{\mathbf G}^{\circ} ([\mathbf L,\mathbf L] ) )^F$.  Hence,  it suffices to prove that  $p$  does not divide   $|M:  (C_{\mathbf G}^{\circ} ([\mathbf L, \mathbf L] ) )^FA|  $.  Since $C _{\mathbf G}^{\circ} ([\mathbf L, \mathbf L] ) $ and  $[\mathbf L, \mathbf L]  $     commute element wise  and  $C _{\mathbf G}^{\circ} ([\mathbf L, \mathbf L] ) \cap  [\mathbf L, \mathbf L]   \leq  Z([\mathbf L, \mathbf L])  $, a standard application of the Lang--Steinberg theorem, see \cite[3.10]{dm},  gives 
$$|
(C_{\mathbf G}^{\circ} ([\mathbf L, \mathbf L] )[\mathbf L,\mathbf L])| ^F  =  |C_{\mathbf G}^{\circ} ([\mathbf L, \mathbf L] )|^F  |[\mathbf L, \mathbf L] |^F  =  |C_{\mathbf G}^{\circ} ([\mathbf L, \mathbf L] )|^F |A|.   $$
On the other hand,   $C_{\mathbf G}^{\circ} ([\mathbf L, \mathbf L] )^ F \cap    A \leq    Z([\mathbf L, \mathbf L]^F)  = Z([\mathbf L, \mathbf L]) ^ F $  is a  $p' $-group  by \cite[Proposition 4]{newf}.    Hence,   the $p$-part  of    $      |C_{\mathbf G}^{\circ} ([\mathbf L, \mathbf L] )^F|   |A|     =  \frac  {   |C_{\mathbf G}^{\circ} \mathbf L, \mathbf L] )^F|  |A| } { |C_{\mathbf G}^{\circ} ([\mathbf L, \mathbf L] ^ F \cap    A |   }   $  equals   the $p$-part of    $ |C_{\mathbf G}^{\circ} ([\mathbf L, \mathbf L] )^F|  |A| $. 
Thus, by the above displayed equation,  it suffices to prove that  $p$  does not divide   $|M:   (C_{\mathbf G}^{\circ} ([\mathbf L, \mathbf L] ) [\mathbf L, \mathbf L] ) ^F|$.   This follows from \cite[Proposition 6]{newf} 
\hfill$\square$\\

\textit{Proof of Theorem \ref{kombucha}}. Let $H:=N_{\mathbf N^F}([\mathbf L,\mathbf L], \res^{\mathbf L^F}_{[\mathbf L,\mathbf L]^F} \zeta, d)$, where $d$ is the unique block of $kM$, where $M:=N_{\mathbf G^F}([\mathbf L,\mathbf L], \res^{\mathbf L^F}_{[\mathbf L,\mathbf L]^F} \zeta)$, with $(P,e_P)$ a $d$-Brauer pair. By Proposition \ref{potas}, we have $\mathcal F_{(P,e_P)}(\mathbf N^F, \mathbf G^F, b) \cong \mathcal F_{(P,e_P)}(H,M,d)$.\\
Next we apply the Generalised Second Fong Reduction to $\mathcal F_{(P,e_P)}(H,M,d)$. For this, let $A:=[\mathbf L,\mathbf L]^F$. In particular, we have normal inclusions $A \leq M \leq H$ with $A \unlhd H$ and $|H:M|_p=1$. Let $c$ be the block of $A$ containing $\res^{\mathbf L^F}_{[\mathbf L,\mathbf L]^F} \zeta$. Apply \cite[Theorem 22.9(ii)]{caen} to $\mathbf L$ to see that $\res^{\mathbf L^F}_{[\mathbf L,\mathbf L]^F} \zeta$ is a defect zero character and thus $c$ is of defect zero. Clearly, $c$ is $M$-stable, so we can apply Lemma \ref{coverd} with $Z=Z(\mathbf L) _{p}^F$ and $e=e_Z$ to see that $c$ is covered by $d$. Thus, we can apply Theorem \ref{c} to this situation and get a fusion system $\mathcal F_{(P,e_P')}(\widetilde H,\widetilde M,\tilde d)$ such that $\widetilde H$ is a $p'$-central extension of $H/A$, $\widetilde M$ is the full inverse image of $M/A$ in $\widetilde H$ and $\tilde d$ is an $\widetilde H$-stable block of $k \widetilde M$ with maximal $(\tilde d, \widetilde H)$-Brauer pair $(P,e_P')$.\\
The final step is now to prove that $\mathcal F_{(P,e_P')}(\widetilde H,\widetilde M,\tilde d)$ is non-exotic. Let $X \leq P$ be $\mathcal F_{(P,e_P')}(\widetilde H,\widetilde M,\tilde d)$-centric, then by Lemma \ref{2020} we have that the normaliser of $(X,e_X)$ is the normaliser of $X$. Furthermore, $P \in \syl_p(\widetilde H)$ by Lemma \ref{2000}, which together with the observation about the normalisers we just mentioned implies $\mathcal F_{(P,e_P')}(\widetilde H,\widetilde M,\tilde d)=\mathcal F_{P}(\widetilde H)$, which means that the generalised block fusion system is non-exotic. \hfill$\square$

\subsection{Remark} Comparison with Theorem \ref{virginia}: We discuss how the proof of Theorem \ref{kombucha} differs from the proof Theorem \ref{virginia}. Firstly, note that the proof of Theorem \ref{kombucha} does not use Theorem \ref{virginia} directly, thus it gives an alternative, more general proof for this theorem. The proof of Theorem \ref{virginia} uses a concept called control subgroups of fusion systems. These groups are defined in \cite{thevenaz} and are calculated for fusion systems of unipotent blocks in \cite{newf}. Let $(\mathbf L, \zeta)$ be a unipotent $e$-cuspidal pair of $\mathbf G$ defining a unipotent block $b$ of $\mathbf G^F$. Cabanes and Enguehard show that $H=(N_V^F)SC_{\mathbf G}^\circ([\mathbf L,\mathbf L])^F$ is a control subgroup of the fusion system of $b$. We discuss the factors that make up this group. Let $\bf T$ be a maximally split torus of $[\mathbf L, \mathbf L]$ $C_{\mathbf G}^\circ([\mathbf L,\mathbf L])$ and $\mathbf K=[\mathbf L, \mathbf L]$ $C_{\mathbf G}([\mathbf L,\mathbf L])$. The group $V$ is the stabiliser of an $F$-stable basis of the root system of $\mathbf K$ with respect to $\mathbf T$ in $N_M(\mathbf T)/\mathbf T^F$, where $M=N_{\mathbf G^F}([\mathbf L,\mathbf L],\res^{\mathbf L^F}_{[\mathbf L,\mathbf L]^F} \zeta)$. Now $N$ is a subgroup of $N_{\mathbf G}(\mathbf T)$ such that $N \cap \mathbf T = \mathbf T_{\phi_2}$ and $N T = N_{\mathbf G}(\mathbf T)$. Here $\mathbf T_{\phi_2}$ denotes the $\phi_2$-part of $\bf T$ in the sense of the notation given after Theorem \ref{naomiwatts}. The existence of the group $N$ is proven in \cite[Proposition 4.2]{cabanes}, which relies on \cite[4.4]{tits}. The group $N_V$ is defined to be the inverse image of $V$ in $N$ under the quotient map by $\mathbf T$. Finally, $S=S'[S',V]$ for some $S' \leq \mathbf T^F_{p'}$ which is characterised by the equation $\mathbf T^F=S'(\mathbf T \cap [\mathbf L,\mathbf L]^FC_{\mathbf G}([\mathbf L,\mathbf L])^F)$.\\

Proving Theorem \ref{kombucha} with control subgroups does not seem to be possible currently. To prove that the control subgroup is invariant under extending by $\bf N$ as in the assumptions of the theorem, we would need to know whether the group $N$ is unique up to conjugation by an element of $N_{\mathbf G}(\mathbf T)$. This boils down to the uniqueness up to conjugation of the Tits extension, see \cite[4.4]{tits}. Alternatively the existence of such an extension for disconnected reductive groups would also imply the existence of such a group $N$ in the situation of the theorem. However, both questions are currently not answered and thus it is unclear how to prove the theorem in the same way as in \cite{oh}. The proof of Theorem \ref{kombucha} circumvents using control subgroups and thus Tits’ theorem for splitting. Here we construct the realising group indirectly, not as a subgroup of $\mathbf G^F$.

\section*{Acknowledgements}

I thank Radha Kessar and Gunter Malle. I also thank the referees for their careful reading of my manuscript and their many valuable comments.

\end{document}